\documentclass[12pt]{article}
\usepackage{latexsym,amsfonts,amssymb,amsmath}
\usepackage{color}
\usepackage{hyperref}
\usepackage{cleveref}

\usepackage{colordvi,multicol}

\newtheorem{thm}{Theorem}[section]

\newtheorem{lem}[thm]{Lemma}
\newtheorem{pro}[thm]{Proposition}
\newtheorem{defi}[thm]{Definition}

\date{}
\begin{document}

\title{\bf A study on the  Weibull and Pareto  distributions motivated by Chv\'{a}tal's theorem}
 \author{Cheng Li$^1$, Ze-Chun Hu$^{1}$, Qian-Qian Zhou$^{2,}$\footnote{Corresponding author: qianqzhou@yeah.net}\\ \\
   {\small $^1$ College of Mathematics, Sichuan  University,
 Chengdu 610065, China}\\
 {\small $^2$ School of Science, Nanjing University of Posts and Telecommunications, Nanjing  210023, China}}
\maketitle

\begin{abstract}

  Let $B(n,p)$ denote a binomial random variable with parameters $n$ and $p$.  Chv\'{a}tal's theorem  says that for any fixed $n\geq 2$, as $m$ ranges over $\{0,\ldots,n\}$, the probability $q_m:=P(B(n,m/n)\leq m)$ is the smallest when  $m$ is closest to $\frac{2n}{3}$.  Motivated by this theorem,  we consider the minimum value problem on the probability that a random variable is at most its expectation, when its distribution is the Weibull distribution or the Pareto distribution in this note.
\end{abstract}



\section{Introduction}

Let $B(n,p)$ denote a binomial random variable with parameters $n$ and $p$. Janson in \cite{Ja21} introduced the following conjecture suggested by Va\v{s}k Chv\'{a}tal.

{\bf Conjecture 1} (Chv\'{a}tal). For any fixed $n\geq 2$, as $m$ ranges over $\{0,\ldots,n\}$, the probability $q_m:=P(B(n,m/n)\leq m)$ is the smallest when
 $m$ is closest to $\frac{2n}{3}$.

Conjecture 1 has significant   applications in machine learning, such as  the analysis of generalized boundaries of relative deviation bounds and unbounded loss functions (\cite{Do18} and \cite{GM14}).
As to the probability of a binomial random variable exceeding its expectation, we refer to Doerr \cite{Do18}, Greenberg  and Mohri \cite{GM14}, Pelekis and Ramon \cite{PR16}.
 Janson \cite{Ja21} proved that  Conjecture 1 holds for large $n$. Barabesi et al. \cite{BPR21} and Sun \cite{Su21} gave an affirmative answer to Conjecture 1 for general $n\geq 2$. Hereafter, we call Conjecture 1 by Chv\'{a}tal's theorem.

Motivated by Chv\'{a}tal's theorem,  Li et al. \cite{LXH23}  considered the minimum value problem on the probability that a random variable is not more than its expectation, when its distribution is the Poisson distribution, the geometric distribution or the Pascal distribution. Sun et al. \cite{SHS23} investigated  the corresponding minimum value problem for the Gamma distribution among other things.  In this note, we consider  the minimum value problem for the Weibull distribution and the Pareto distribution in Sections 2 and 3, respectively.

\section{The Weibull  distribution}\setcounter{equation}{0}

Let $X$ be a  Weibull random variable with parameters $\alpha$ and $\theta\,(\alpha > 0,\theta > 0)$ and the density function
$$f(x)=\frac{\alpha }{\theta }x^{\alpha -1}e^{-\frac{x^{\alpha}}{\theta }},\, x>0.
$$
We know that its expectation $EX=\theta ^{\frac{1}{\alpha }}\Gamma \left ( \frac{1}{\alpha } +1\right ),$  where $\Gamma \left ( \frac{1}{\alpha } +1\right )$  is the Gamma function, i.e., $\Gamma \left ( \frac{1}{\alpha } +1\right )=\int_0^{\infty}u^{\frac{1}{\alpha}}e^{-u}du$. For any given real number $\kappa>0,$   we have
$$P\left ( X\leq \kappa  EX \right )=\int_{0}^{\kappa \theta ^{\frac{1}{\alpha }}\Gamma \left ( \frac{1}{\alpha } +1\right )}\frac{\alpha }{\theta }t^{\alpha -1}e^{-\frac{t^{\alpha}}{\theta }}dt.
$$
By taking the change of variable $t= \left ( \theta x \right )^{\frac{1}{\alpha }}$, we get
\begin{eqnarray*}
P\left (X\leq \kappa  EX \right )&=& \int_{0}^{\left (\kappa
 \Gamma \left ( \frac{1}{\alpha } +1\right )\right)^{\alpha }}\frac{\alpha }{\theta }\left ( \theta x \right )^{\frac{\alpha-1}{\alpha }}e^{-x} \frac{\theta ^{\frac{1}{\alpha }}}{\alpha }x^{\frac{1}{\alpha }-1}dx \\
&=&\int_{0}^{\left (\kappa \Gamma \left ( \frac{1}{\alpha } +1\right )\right)^{\alpha }}e^{-x}dx\\
& =&1-e^{-\left (\kappa \Gamma \left ( \frac{1}{\alpha } +1\right )\right)^{\alpha }},
\end{eqnarray*}
which shows that $P\left ( X\leq \kappa EX \right )$ is independent of $\theta$.

Define a function
\begin{eqnarray}\label{3.1}
g_{\kappa}( \alpha ):=1-e^{-\left (\kappa \Gamma \left ( \frac{1}{\alpha } +1\right )\right)^{\alpha }},\quad \alpha > 0.
\end{eqnarray}

The main result of this section is

\begin{pro}\label{pro-2.1}
(i) If $\kappa \leq 1$, then
\begin{eqnarray*}
\inf\limits_{\alpha \in \left ( 0,+\infty  \right ) }g_{\kappa }\left ( \alpha  \right )=\lim_{\alpha\to +\infty}g_{\kappa}(\alpha)=
\left\{
\begin{array}{cl}
0,& \kappa < 1,\\
1-e^{-e^{-\gamma }},& \kappa = 1,
\end{array}
\right.
\end{eqnarray*}
where $\gamma$  is the Euler's constant, i.e., $\gamma=\sum_{n=1}^{\infty }\left [ \frac{1}{n} -\ln\left ( 1+\frac{1}{n} \right )\right ]$.\\
(ii) If $\kappa>1$, then
$$
\min_{\alpha \in \left ( 0,+\infty  \right ) }g_{\kappa }\left ( \alpha  \right )=g_{\kappa }\left ( \alpha _{0}\left ( \kappa  \right )  \right ),
$$
where $\alpha _{0}\left ( \kappa  \right )=\frac{1}{x_{0}\left ( \kappa  \right ) -1}$, and $x_{0}\left ( \kappa  \right )$  is the unique null point of function $\varphi _{\kappa }\left ( x \right ):= \left ( x-1 \right )\psi \left ( x \right )-\ln\left ( \kappa \Gamma \left ( x \right ) \right )$ on  $\left ( 1,+\infty  \right ),$  where $\psi(x)$ is the digamma function (see Definition 2.3 below).
\end{pro}

Note that $\left (\kappa \Gamma \left ( \frac{1}{\alpha } +1\right )\right)^{\alpha }=e^{\alpha \ln\left (\kappa \Gamma \left ( \frac{1}{\alpha }+1 \right )\right)}.$  Let $x=\frac{1}{\alpha }+1,$   and define function
\begin{eqnarray}\label{3.2}
h_{\kappa}(x):=\frac{\ln (\kappa \Gamma(x))}{x-1},\quad x > 1.
\end{eqnarray}
Then
\begin{eqnarray}\label{3.3}
g_{\kappa }\left ( \alpha  \right )=1-e^{-e^{h_{\kappa}\left ( x  \right )}},
\end{eqnarray}
 and in order to finish the proof of Proposition \ref{pro-2.1}, it is enough to prove the following lemma.

\begin{lem} \label{lem-2.2}
(i) If $\kappa \leq 1,$  then
\begin{eqnarray*}
\inf\limits_{x\in \left ( 1,+\infty  \right ) }h_{\kappa}\left ( x  \right )=\lim_{x\to 1^+}h_{\kappa}(x)=
\left\{
\begin{array}{ll}
-\infty,& \kappa < 1,\\
-\gamma,& \kappa = 1,
\end{array}
\right.
\end{eqnarray*}
where $\gamma$ is the Euler's constant.\\
(ii) If $\kappa>1$, then
$$\min\limits_{x\in \left ( 1,+\infty  \right ) }h_{\kappa}\left ( x  \right )=h_{\kappa}\left ( x_{0}\left ( \kappa  \right ) \right ),$$
where $x_{0}\left ( \kappa  \right )$  is the unique null point of function $\varphi _{\kappa }\left ( x \right ):= \left ( x-1 \right )\psi \left ( x \right )-\ln\left ( \kappa \Gamma \left ( x \right ) \right )$ on $\left ( 1,+\infty  \right ),$  where $\psi(x)$ is the digamma function.
\end{lem}

Before giving the proof of Lemma \ref{lem-2.2}, we need some preliminaries on ploygamma function. 

\begin{defi}\label{defi-3.2} (\cite[1.16]{Er55}) Let $m$ be any nonnegative integers. $m$-order ploygamma function $\psi^{(m)}$ is defined by
\begin{eqnarray*}
\psi ^{(m)}( z) := \frac{\mathrm{d}^{m} }{\mathrm{d} z^{m}}\psi (z)=\frac{\mathrm{d}^{m+1} }{\mathrm{d} z^{m+1}}\ln\Gamma (z), Re z>0.
\end{eqnarray*}
When $m=0$, $\psi(z):=\psi ^{( 0 )}( z)=\frac{\mathrm{d} }{\mathrm{d} z}\ln\Gamma( z)=\frac{{\Gamma}'(z)}{\Gamma(z)}$ is called digamma function.
\end{defi}

By \cite[1.7(3)]{Er55} and \cite[1.9(10)]{Er55}, we know that
\begin{eqnarray}
\psi(z)&=& -\gamma -\frac{1}{z}+\sum_{n=1}^{\infty }\frac{z}{n(z+n)}\nonumber\\
& =&-\gamma +( z-1)\sum_{n=0}^{\infty }\frac{1}{[( n+1)( z+n )]},\label{3.2-a}\\
\psi ^{( 1 )}( z) &=&{\psi }'( z)=\sum_{n=0}^{\infty }\frac{1}{( z+n)^{2}}.\label{3.2-b}
\end{eqnarray}

\noindent {\bf Proof of Lemma \ref{lem-2.2}.} By \eqref{3.2} and Definition \ref{defi-3.2}, we have
\begin{eqnarray}\label{3.7}
{h}'_{\kappa }\left ( x \right )=\frac{\left ( x-1 \right )\psi \left ( x \right )-\ln\left ( \kappa \Gamma \left ( x \right ) \right )}{\left ( x-1 \right )^{2}}=\frac{\varphi _{\kappa }\left ( x \right )}{\left ( x-1 \right )^{2}},\quad x> 1.
\end{eqnarray}
By \eqref{3.2-b}, we get
\begin{eqnarray*}
{\varphi}'_{\kappa }\left ( x \right )=\left ( x-1 \right )\psi ^{\left ( 1 \right )}\left ( x \right )=\left (x-1 \right )\sum_{n=0}^{\infty }\frac{1}{\left ( x+n \right )^{2}}>0,\quad \forall x>1.
\end{eqnarray*}
It follows that the function $\varphi _{\kappa } \left ( x \right )$ is strictly increasing on the interval $\left ( 1,+\infty  \right )$.

Thus, if $\kappa \leq 1$, we have
$$\varphi _{\kappa }\left ( x \right )> \varphi _{\kappa }\left ( 1 \right )=-\ln\kappa \geq 0, \quad \forall x >1.$$
Then, by \eqref{3.7}  we get
$${h}'_{\kappa }\left ( x \right )>0,\quad \forall x >1,$$
which implies that the function  $h_{\kappa}\left ( x \right )$ is strictly increasing on $\left ( 1,+\infty  \right ).$   Hence the function $h_{\kappa}\left ( x \right )$ has no  minimum value on $\left ( 1,+\infty  \right )$ and
$$
\inf\limits_{x\in \left ( 1,+\infty  \right ) }h_\kappa\left ( x  \right )=\lim\limits_{x\rightarrow 1^+}h_{\kappa}\left ( x  \right )=\lim\limits_{x\rightarrow 1^+}\frac{\ln\Gamma \left ( x \right )}{x-1}+\lim\limits_{x\rightarrow 1^+}\frac{\ln\kappa }{x-1}.
$$
By the L'Hospital's rule and \eqref{3.2-a}, we have
\begin{eqnarray*}
\lim\limits_{x\rightarrow 1^+}\frac{\ln\Gamma \left ( x \right )}{x-1}& =\lim\limits_{x\rightarrow 1^+}\frac{\Gamma'(x)}{\Gamma(x)}={\Gamma }'\left ( 1 \right )=\psi \left ( 1 \right )=-\gamma.
\end{eqnarray*}
Thus, 
\begin{eqnarray*}
\lim_{x\in (1,+\infty) }h_{\kappa}(x)=
\left\{
\begin{array}{ll}
-\infty,& \kappa < 1,\\
-\gamma,& \kappa = 1.
\end{array}
\right.
\end{eqnarray*}

If $\kappa>1$, then $\varphi _{\kappa }\left ( 1 \right )=-\ln\kappa <0.$  By  \cite[1.18(1)]{Er55} (Stirling formula) and  \cite[1.18(7)]{Er55},
 when $z\rightarrow \infty$ we have
\begin{eqnarray*}
\ln\Gamma \left ( z \right )&= \left ( z-\frac{1}{2} \right )\ln z-z+\frac{\ln\left ( 2\pi  \right )}{2}+o\left ( 1 \right ),\\
\psi \left ( z \right )&= \ln z-\frac{1}{2z}+o\left ( \frac{1}{z} \right ),\quad \left | \arg z \right |< \pi .
\end{eqnarray*}
Then
\begin{eqnarray*}
\lim_{x\rightarrow +\infty}\varphi_{\kappa}(x)&=&\lim_{x\rightarrow +\infty}[(x-1)\psi(x)-\ln \Gamma(x)-\ln \kappa]\\
&=&\lim_{x\rightarrow +\infty}\left[ (x-1)\left(\ln x-\frac{1}{2x}+o\left(\frac{1}{x}\right)\right)\right.\\
&&\quad\quad\quad\left.-\left(\left(x-\frac{1}{2}\right)\ln x -x +\frac{\ln (2\pi)}{2}+o(1)\right)-\ln \kappa\right]\\
&=&\lim_{x\rightarrow +\infty}\left[ x-\frac{1}{2}\ln x +\frac{1}{2x}-\frac{1}{2}-\frac{\ln\left ( 2\pi  \right )}{2}-\ln\kappa +o\left ( 1 \right )   \right]\\
&=&+\infty.
\end{eqnarray*}
Since the function $\varphi _{\kappa } \left ( x \right )$ is continuous, by the zero point theorem, there exists $x_{0}\left ( \kappa  \right ) \in (1,+\infty)$ which depends on  $\kappa $ satisfying  that $$\varphi _{\kappa } \left ( x_{0}\left ( \kappa  \right ) \right )=0.$$
Moreover, combining with  the monotonicity of the function $\varphi _{\kappa } \left ( x \right )$ on the interval  $\left ( 1,+\infty  \right )$, we know that $x_{0}\left ( \kappa  \right )$ is the unique null point of the function $\varphi _{\kappa } \left ( x \right )$  and
\begin{eqnarray*}
\varphi _{\kappa } \left ( x \right )&<0,\quad \forall x\in \left ( 1,x_{0}\left ( \kappa  \right ) \right );\\
\varphi _{\kappa } \left ( x \right )&>0,\quad \forall x\in \left ( x_{0}\left ( \kappa  \right ),+\infty  \right ).
\end{eqnarray*}
Then, by \eqref{3.7}  we get
\begin{eqnarray*}
{h}'_{\kappa }\left ( x \right )&<0,\quad \forall x\in \left ( 1,x_{0}\left ( \kappa  \right ) \right );\\
{h}'_{\kappa }\left ( x \right )&>0,\quad \forall x\in \left ( x_{0}\left ( \kappa  \right ),+\infty  \right ).
\end{eqnarray*}
Thus, the function $h_{\kappa }\left ( x \right )$ is strictly  decreasing on  $\left ( 1,x_{0}\left ( \kappa  \right ) \right )$  and strictly increasing on $\left ( x_{0}\left ( \kappa  \right ),+\infty  \right )$, which implies that
$$h_{\kappa }\left ( x \right )\geq h_{\kappa }\left ( x_{0}\left ( \kappa  \right ) \right ),\quad \forall x> 1.$$
Therefore,
$$\min\limits_{x\in \left ( 1,+\infty  \right ) }h_{\kappa}\left ( x  \right )=h_{\kappa}\left ( x_{0}\left ( \kappa  \right ) \right ).$$
The proof is complete. \hfill\fbox

\section{The Pareto distribution}\setcounter{equation}{0}

Let $X$ be a  Pareto  random variable with parameters $a$ and $\theta\, (a>0,\theta>0)$ and the density function
$$
f(x)=\theta a^{\theta }x^{-(\theta +1)}I_{(a,\infty)}(x).
$$
 When $\theta>1$, the expectation of $X$ is $EX=\frac{\theta a}{\theta -1}$.
  Then, for any given real number $\kappa>0$, we have
\begin{eqnarray*}
P\left ( X\leq \kappa EX \right )&=&\int_{a}^{\frac{\kappa \theta a}{\theta -1}}\theta a^{\theta }t^{-\left ( \theta +1 \right )}dt\\
&=&-\left ( \frac{a}{t} \right )^{\theta }\big|_{a}^{\frac{\kappa \theta a}{\theta -1}}\\
&=&1-\left ( \frac{\theta -1}{\kappa \theta }\right )^{\theta },
\end{eqnarray*}
which shows that $P\left ( X\leq \kappa EX \right )$ is independent of $a$. Note that, in order to make sense of  the above equality, if $\kappa<1$, the parameter $\theta$ should satisfy that $1<\theta\le \frac{1}{1-\kappa};$  and if $\kappa\ge 1$,  the parameter $\theta$ should  satisfy that $\theta>1.$

Define a function
$$
g_{\kappa}(\theta):=  1-\left(\frac{\theta -1}{\kappa \theta }\right )^{\theta },\ 1<\theta \leq \frac{1}{1-\kappa },\kappa <1\ \mbox{or}\ \theta > 1,\kappa \geq 1.
$$

The main result of this section is

\begin{pro}\label{pro-3.1}
(i) If $\kappa<1$, then
$$
\min\limits_{\theta \in \left ( 1,\frac{1}{1-\kappa }  \right ] }g_{\kappa }\left ( \theta  \right )=g_{\kappa }\left ( \frac{1}{1-\kappa}\right )=0.
$$
(ii) If $\kappa=1$, then $\inf\limits_{\theta \in ( 1,+\infty)}g_1(\theta)=\lim\limits_{\theta\to +\infty}g_1(\theta)=1-e^{-1}$. \\
(iii) If $\kappa > 1$, then
$$
\min\limits_{\theta \in \left ( 1,+\infty  \right ) }g_{\kappa }\left ( \theta  \right )=
g_{\kappa}\left ( \theta_{0}\left ( \kappa  \right ) \right ),
$$
where $\theta _{0}\left ( \kappa  \right )=\frac{1}{1-x_{0}\left ( \kappa  \right ) }$, and $x_{0}\left ( \kappa  \right )$ is the unique null point of function  $\varphi _{\kappa }\left ( x \right ):= 1-\frac{1}{x}-\ln\frac{x}{\kappa }$ on the interval $\left ( 0,1  \right )$.
\end{pro}

Note that $\left ( \frac{\theta -1}{\kappa \theta }\right )^{\theta }=e^{\theta \ln\frac{\theta -1}{\kappa \theta }}$. Let $x=1-\frac{1}{\theta }$
and define function
\begin{eqnarray}\label{4.1}
h_{\kappa }\left ( x \right ):=\frac{\ln \frac{x}{\kappa }}{x-1},\ 0<x \leq \kappa ,\kappa < 1,\ \mbox{or}\ 0<x<1,\kappa \geq 1.
\end{eqnarray}
Then
\begin{eqnarray}\label{4.2}
g_{\kappa}\left ( \theta  \right )=1-e^{-h_{\kappa }\left ( x \right )},
\end{eqnarray}
and in order to finish the proof of Proposition \ref{pro-3.1}, it is enough to prove the following lemma.

\begin{lem}\label{lem-3.2}
(i) If $\kappa<1$, then
$$\min\limits_{x\in ( 0,\kappa] }h_{\kappa}(x)=h_{\kappa}(\kappa)=0.$$\\
(ii) If $\kappa=1$, then $\inf\limits_{x\in ( 0,1) }h_1(x)=\lim\limits_{x\to 1^-}h_1(x)=1$.\\
(iii) If $\kappa>1$, then
\begin{eqnarray*}
\min\limits_{x\in ( 0,1) }h_{\kappa}(x)=h_{\kappa}\left ( x_{0}\left ( \kappa  \right ) \right ),
\end{eqnarray*}
where $x_{0}\left ( \kappa  \right )$ is the unique null point of function  $\varphi _{\kappa }\left ( x \right ):= 1-\frac{1}{x}-\ln\frac{x}{\kappa }$ on the interval $\left ( 0,1  \right )$.
\end{lem}

\noindent {\bf Proof.} (i) If $\kappa <1$, by ~\eqref{4.1} we have
\begin{eqnarray}\label{4.3}
{h}'_{\kappa }\left ( x \right )=\frac{1-\frac{1}{x}-\ln\frac{x}{\kappa }}{\left ( x-1 \right )^{2}}=\frac{\varphi _{\kappa }\left ( x \right )}{\left ( x-1 \right )^{2}},\quad 0<x \leq \kappa.
\end{eqnarray}
By the definition of   $\varphi _{\kappa }\left ( x \right )$, we get that
\begin{eqnarray*}\label{4.4}
{\varphi}'_{\kappa }\left ( x \right )=\frac{1-x}{x^{2}}> 0,\quad \forall 0<x \leq \kappa.
\end{eqnarray*}
It follows  that  function $\varphi _{\kappa }\left ( x \right )$ is strictly increasing on  $\left ( 0,\kappa \right]$ and thus
\begin{eqnarray}\label{4.4-1}
\varphi _{\kappa }\left ( x \right ) \leq \varphi _{\kappa }\left ( \kappa \right )=\frac{\kappa -1}{\kappa }< 0,\quad \forall 0<x\leq \kappa.
\end{eqnarray}
Then, by \eqref{4.3}  and \eqref{4.4-1} we get
$${h}'_{\kappa }\left ( x \right )<0,\quad \forall 0<x\leq \kappa,$$
which implies that  function $h_{\kappa}\left ( x \right )$ is strictly decreasing on  $\left ( 0,\kappa \right]$. Thus
$$\underset{x\in \left ( 0,\kappa \right]}{\min}h_{\kappa}\left ( x  \right )=h_{\kappa}\left ( \kappa \right )=0.$$

 If $\kappa \geq 1$, by~\eqref{4.1} and the definition of  $\varphi_{\kappa}(x)$ again, we also have
\begin{eqnarray}\label{4.5}
{h}'_{\kappa }\left ( x \right )=\frac{\varphi _{\kappa }\left ( x \right )}{\left ( x-1 \right )^{2}},\quad 0< x <1,
\end{eqnarray}
and
\begin{eqnarray*}\label{4.6}
{\varphi}'_{\kappa }\left ( x \right )=\frac{1-x}{x^{2}}> 0,\quad \forall 0 <x <1.
\end{eqnarray*}
It follows that function $\varphi _{\kappa }\left ( x \right )$ is strictly increasing on  $\left ( 0,1 \right)$.

(ii) If $\kappa = 1$, then
$$\varphi _{\kappa }\left ( x \right ) < \varphi _{\kappa }\left ( 1 \right )=-\ln \kappa= 0,\quad \forall 0<x<1.$$
By \eqref{4.5}, we get that
 $${h}'_{\kappa }\left ( x \right )<0,\quad \forall 0<x<1,$$
which implies  that function $h_{\kappa}\left ( x \right )$ is strictly decreasing on   $\left ( 0,1 \right)$. Thus,
$$\inf\limits_{x\in \left ( 0,1 \right) }h_{\kappa}\left ( x  \right )=\lim\limits_{x\rightarrow 1^-}h_{\kappa}\left ( x  \right )=\lim\limits_{x\rightarrow 1^-}\frac{\ln x}{x-1}=1.$$

(iii)  If $\kappa >1$, then $\varphi _{\kappa }\left ( 1 \right )=\ln \kappa>0$. Moreover,
\begin{eqnarray*}
\lim\limits_{x\rightarrow 0^+ }\varphi _{\kappa } \left ( x \right )&=\lim\limits_{x\rightarrow 0^+ }\left(1-\frac{1}{x}-\ln x+\ln\kappa \right)\\
&=\lim\limits_{x\rightarrow 0^+ }\left ( 1+\ln\kappa -\frac{x\ln x+1}{x} \right )\\
&=-\infty.
\end{eqnarray*}
Since the function $\varphi _{\kappa } \left ( x \right )$ is continuous on $(0,1)$, by  the zero point theorem, there exists  $x_{0}\left ( \kappa  \right ) \in (0,1)$ depending on parameter $\kappa $   fulfills that $$\varphi _{\kappa } \left ( x_{0}\left ( \kappa  \right ) \right )=0.$$
By the monotonicity of function $\varphi _{\kappa } \left ( x \right )$ on $\left ( 0,1 \right)$, we know that $x_{0}\left ( \kappa  \right )$ is the unique  null point of  $\varphi _{\kappa } \left ( x \right )$ and
\begin{eqnarray*}
\varphi _{\kappa } \left ( x \right )&<0,\quad \forall x\in \left ( 0,x_{0}\left ( \kappa  \right ) \right );\\
\varphi _{\kappa } \left ( x \right )&>0,\quad \forall x\in \left ( x_{0}\left ( \kappa  \right ),1  \right ).
\end{eqnarray*}
Then, by \eqref{4.5} we have
\begin{eqnarray*}
{h}'_{\kappa }\left ( x \right )&<0,\quad \forall x\in \left ( 0,x_{0}\left ( \kappa  \right ) \right );\\
{h}'_{\kappa }\left ( x \right )&>0,\quad \forall x\in \left ( x_{0}\left ( \kappa  \right ),1  \right ).
\end{eqnarray*}
Therefore, the function $h_{\kappa }\left ( x \right )$ is strictly decreasing on  $\left ( 0,x_{0}\left ( \kappa  \right ) \right )$ and is strictly increasing on  $\left ( x_{0}\left ( \kappa  \right ),1  \right )$.
Thus
$$\underset{x\in \left ( 0,1\right )}{\min}h_{\kappa}\left ( x  \right )=h_{\kappa}\left ( x_{0}\left ( \kappa  \right ) \right ).$$
The proof is complete. \hfill\fbox

\bigskip

\noindent {\bf\large Acknowledgments}\quad  This work was supported by the National Natural Science Foundation of China (12171335), the Science Development Project of Sichuan University (2020SCUNL201) and  the Scientific Foundation of Nanjing University of Posts and Telecommunications (NY221026).

\end{document}